\documentclass[conference]{IEEEtran}
\IEEEoverridecommandlockouts
\usepackage{cite}
\usepackage{amsmath,amssymb,amsfonts}
\usepackage{algorithmic}
\usepackage{cite}
\ifCLASSINFOpdf
   \usepackage[pdftex]{graphicx}
\else
   \usepackage[dvips]{graphicx}
\fi

\usepackage{booktabs} 
\usepackage{array} 
\usepackage{tabularx}
\usepackage{textcomp}
\usepackage{booktabs}
\usepackage{siunitx}
\newcolumntype{C}[1]{>{\centering\arraybackslash}p{#1}}
\usepackage{array} 
\usepackage{eurosym} 
\usepackage{mathtools}
\usepackage{subcaption}
\usepackage{hyperref}
\usepackage{amsmath}
\usepackage{amssymb}
\usepackage{accents}

\usepackage{xcolor}
\def\BibTeX{{\rm B\kern-.05em{\sc i\kern-.025em b}\kern-.08em
    T\kern-.1667em\lower.7ex\hbox{E}\kern-.125emX}}
\begin{document}


\title{\fontsize{21}{23}\selectfont
Averaging favors MPC: How typical evaluation setups\\
overstate MPC performance for residential battery scheduling
}

\author{\fontsize{12}{14}\selectfont Janik Pinter\IEEEauthorrefmark{1}, Maximilian Beichter\IEEEauthorrefmark{1}, Ralf Mikut\IEEEauthorrefmark{1}, Frederik Zahn\IEEEauthorrefmark{1}, and Veit Hagenmeyer\IEEEauthorrefmark{1}
\\ \textit{\fontsize{10}{12}\selectfont Institute for Automation and Applied Informatics,}
\\ \textit{\fontsize{10}{12}\selectfont Karlsruhe Institute of Technology (KIT), Germany}
\\ \IEEEauthorrefmark{1}{\tt\small \{firstname.lastname\}@kit.edu}
}

\maketitle

\begin{abstract}
Residential prosumers with PV-battery systems increasingly manage their electricity exchange with the power grid to minimize costs. 
This study investigates the performance of Model Predictive Control (MPC) and Rule-Based Control (RBC) under 15/30/60 minute averaging commonly used in research, when Net Billing and battery degradation are considered. 
We simulate five consecutive months for 15 buildings in northern Germany, generating costs at up to 1-minute resolution while scheduling at 15/30/60 minutes. We find that time-averaged evaluations make MPC look consistently better than RBC, yet when costs are recomputed at minute-level ground-truth, the reported advantage shrinks by 69\% on average for hourly schedulers. For individual buildings, the finer evaluation can reverse conclusions, and simple RBC can achieve lower total costs than an MPC with perfect foresight. These findings caution against drawing conclusions from coarse averages and show how a fair assessment of battery scheduling approaches can be obtained.
\end{abstract}




\section{Introduction}

With the increasing roll-out of residential PV-Battery systems, typical consumers are increasingly turning into prosumers that actively manage their electricity exchange with the power grid.  
Regulators aim to influence the prosumers energy management to improve grid friendliness and with that reduce the necessity of expensive grid expansions while keeping prosumer economics transparent and fair.
For that, their key mechanisms are settlement rules (how imports/exports are accounted for) and tariff structures (how they are priced), which together determine economic incentives and thereby heavily influence the prosumers' energy management strategy. 
Two of the most common approaches for the prosumers' battery control are Rule-Based Control (RBC) and Model Predictive Control (MPC), which is why many studies compare their performance. 
In this paper, we show how such a comparison is often systematically biased to overestimate the advantages of MPC over RBC. 
The bias stems from averaging and misaligned timescales between evaluation practices and real-world settlement rules.
But in order to fully understand when such a bias occurs, it is important to understand the nuances of different metering and billing settlement rules.

\subsection{Metering and Billing Settlement Rules}
Metering and billing settlement rules define the frequency at which the amount of exchanged energy with the grid is measured (e.g., instantaneous, hourly, monthly, yearly) and how those measurements are billed \cite{ziras_effect_2021}.
The main forms are: 
\begin{itemize}
    \item \textbf{Buy All - Sell All}: All generated electricity is directly fed into the grid (typically under a fixed export price for feed-in). All consumed energy is bought from the grid. Self-consumption is not possible.
    \item \textbf{Net Metering}: Imports $a$ and exports $b$ within a fixed time period called Netting Interval (NI) cancel each other out. Only the net exchange $|a-b|$ over the NI is billed.
    \item \textbf{Net Billing}: Imports $a$ and exports $b$ are measured separately and $|a|$ is billed and $|b|$ credited. Self-consumption is possible. Net Billing can be seen as a special form of Net Metering with a very short (instantaneous) NI \cite{ziras_effect_2021}.
\end{itemize}
In this work, we focus on the differences between Net Metering and Net Billing schemes. As of today, there is an ongoing trend towards Net Billing schemes. For example, Denmark introduced Net Metering with an NI of one year in 2010 \cite{gram2020danish}; because settlement was based on the yearly net $|a-b|$, recorded grid interaction---and thus taxable energy volumes---appeared much smaller. Such a system disproportionally favors prosumers over consumers, which is why, since 2012, new residential PV systems in Denmark are subject to Net Billing instead. 
Furthermore, the EU Electricity Directive \cite{eu-directive} started demanding states to charge prosumers for imports and exports separately since 2019. However due to slow implementation across Europe, NI choices still range from very long banking periods (up to 25 years in Greece), through sub-hourly windows (15 minutes in Portugal) to Net Billing (Italy, Poland, Germany, Denmark, Austria, etc.).
See \cite{varga2024changing} for an overview of selected European states and their current settlement rules.

Because the choice of the settlement rule is fundamental for the electricity bill, cost minimization studies should be aware and should make the reader aware of the specific settlement rule they are investigating. Yet, studies rarely specify the settlement rule on which they base their investigation. Although Net Billing is most common, evaluations often average power over 15, 30 or 60 minutes. This effectively sets the settlement rule to Net Metering with a NI equal to the averaging window, which can misalign research with real-world settlement.




\subsection{Related Works}
\label{sec:related-works}

Many studies compare MPC and RBC for residential batteries and typically report notable MPC savings, often assuming perfect forecasts. These results usually arise when evaluation relies on power averaged over 15, 30, or 60 minute windows aligned with the controller's scheduling step \cite{lee_optimal_2018, banfield_comparison_2020, salpakari_optimal_2016, becchi_computationally_2024, benalcazar_transitioning_2024}.
This practice is convenient since it reduces implementation efforts and requires coarser temporal data resolution.
However, computing costs from interval averages implicitly assumes Net Metering with a NI equal to the scheduler's step length, i.e., most often 15/30/60-min. 
Under Net Billing, where imports and exports are priced separately, such averaging masks intra-interval fluctuations that add costs for the common case that the export price is lower than the import price. For example, \cite{banfield_comparison_2020} compare RBC and MPC on 30-min averages considering battery degradation and conclude that MPC with perfect forecasts reduces total costs by 12.9 \% relative to RBC. By only using averages, their cost calculation translates into implicitly assuming that either the load and PV will be constant over each 30-min interval (which does not happen in real-world), or that the evaluation is based on a Net Metering rule with a NI of 30 minutes.

Beyond single-layer comparisons, several works adopt hierarchical control with a slow scheduling layer and a faster execution layer. The fast layer often focuses on tracking a committed dispatch plan (minimizing plan-tracking error) rather than minimizing monetary costs under Net Billing \cite{pasqui_self-dispatching_2025, appino_use_2018, werling_towards_2022}. The authors of \cite{bird_lifetime_2025} implicitly acknowledge the importance of fluctuations by combining a 30-min MPC scheduler with 5-min application. This choice allows for better approximation of costs under a Net Billing scheme. In their setup, the battery power is held constant over an MPC step, so that the 5-min fluctuations are compensated for by the grid. This can be favorable for battery lifetime but may increase the electricity bill significantly due to alternating import/export swings, as our results show later.

Overall, related works commonly resort to averaged evaluation and thus conceal three issues that matter under the most prevalent Net Billing: (i) the implicit assumption about a NI equal to the schedulers step, (ii) the additional cost impact of import/export crossings within an interval, and (iii) the additional cost impact of battery degradation in faster execution layers. 
Whether and how much these factors impact the performances of MPC and RBC differently is not quantified in existing studies. Addressing this gap is the aim of the present work.

\subsection{Contribution}
We investigate the performance of RBC and MPC for the operation of residential PV-Battery systems with respect to cost minimization on different timescales. 
The comparison is based on simulations with real-world measurements from 15 buildings over five consecutive months. 
The costs we evaluate include the electricity bill and the cost of battery degradation. 

These simulations address the mismatch between common averaged evaluations and the dominant settlement practice under Net Billing. The goal is to isolate and measure timescale effects stemming from averages and provide insights into how reliable conclusions drawn from average power evaluations actually are. 
Our key contributions include the following:
\begin{itemize}
    \item Identify Net Metering versus Net Billing mismatch: We create awareness of how averaging over 15, 30, 60 minute windows implicitly assumes Net Metering settlements when insights are intended for Net Billing.

    \item Quantify cost underestimation under coarse time-averages: 
    We quantify how averaging underestimates costs and provide estimates that let researchers roughly approximate Net Billing costs from coarse results.

    \item Expose common evaluation bias: We demonstrate that averaging over 15, 30, or 60 minutes systematically overstates MPC performance relative to RBC under Net Billing, which can lead to distorted or wrong conclusions.  
\end{itemize}

The remainder of the paper is organized as follows: \autoref{sec:prob-formulation} introduces the problem. In \autoref{sec:RBC-MPC}, the investigated RBC and MPC models are presented. The experimental setup, selected battery specifications, pricing tariffs, and data are part of \autoref{sec:exp-setup}. \autoref{sec:results} contains the results of the experiments. \autoref{sec:discussion} discusses the obtained insights. \autoref{sec:conclusio} concludes the paper.

\section{Problem Formulation}
\label{sec:prob-formulation}

\begin{figure}[ht]
    \centering
    \includegraphics[width=0.56\linewidth]{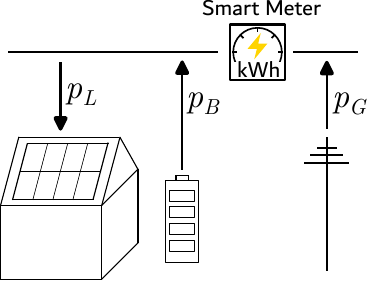}
    \caption{General setting. The residential net-load $p_L$, the battery power $p_B$ and the grid power $p_G$ are in balance.}
    \label{fig:power-balance}
\end{figure}

We investigate the combination of single residential buildings equipped with PV panels and a battery system. We consider the point-of-view of the residents and try to leverage the battery system to minimize our costs, considering the electricity bill and battery degradation. The general structure of this setup is sketched in \autoref{fig:power-balance}. Here, the net-load $p_L(k)$ is shared between the controllable battery power $p_B(k)$ and the grid power $p_G(k)$ such that
\begin{equation}
\label{eq:power-balance}
    p_L(k) = p_B(k) + p_G(k).
\end{equation}
We assume that the battery power $p_B(k)$ is set by a controller. 
The net-load $p_L(k)$ is given from time series measurement data. 
The grid power $p_G(k)$ thus results from the difference between $p_L(k)$ and $p_B(k)$.

\begin{figure}[ht]
    \centering
    \includegraphics[width=0.95\linewidth]{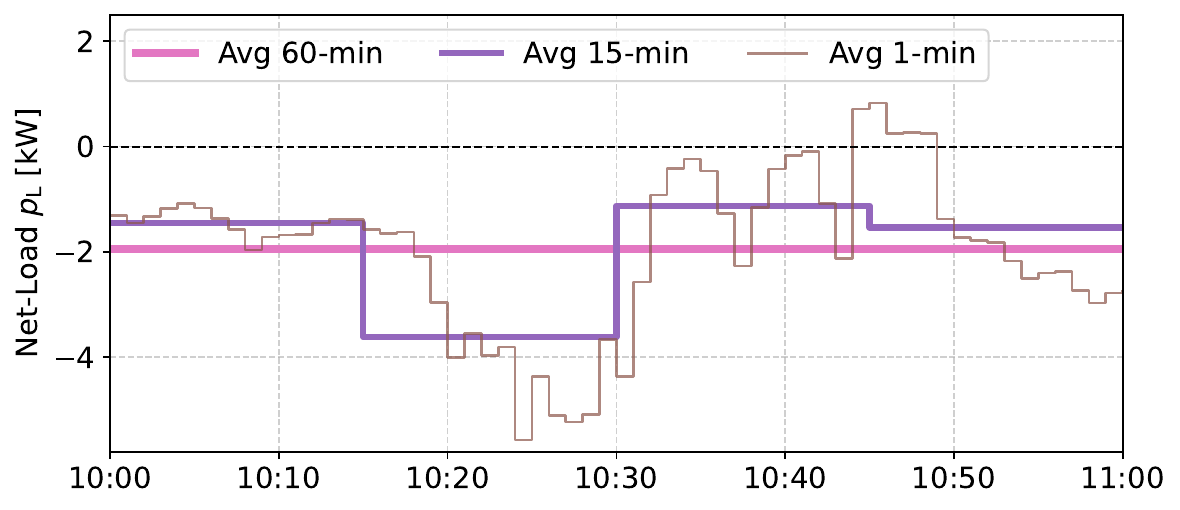}
    \caption{Net-load $p_L$ averaged over three resolutions. High-resolution data can better approximate costs under Net Billing.}
    \label{fig:toy-example}
\end{figure}

In this study, all powers are averaged over predefined time intervals, with $k$ denoting the index of the interval. Because this paper revolves around averaging effects, \autoref{fig:toy-example} illustrates different averaging choices for the same underlying net-load. This illustrates how different averaging choices can cause different battery decisions and hence different costs. The smaller the averaging window, the better the costs can be approximated under Net Billing. 

We formalize this with two timescales: $\Delta gt$ for the ground-truth data resolution and $\Delta NI$ for the Netting Interval used for settlement. Since we focus on the most prevalent Net Billing settlement, we set $\Delta NI = \Delta gt$. As a result, the finest data resolution selected for an experiment directly sets the Netting Interval in the cost calculation. The smaller $\Delta gt$ (and thus $\Delta NI$), the better the cost approximation for a Net Billing settlement rule. We next describe the control schemes that determine the controllable battery power $p_B(k)$.

\section{Battery Control Schemes}
\label{sec:RBC-MPC}

This section outlines the battery control approaches investigated in this study. We first introduce a purely Rule-Based Control (RBC) approach. We then present a two-stage Model Predictive Control (MPC) approach that combines slower forecast-based planning with faster rule-based execution.

\subsection{Rule-Based Control}
Typical RBC is a fast-acting reactive control strategy that controls battery power without relying on forecasts or price information.
We implement the most commonly used rules (see e.g. \cite{banfield_comparison_2020, muller_trade-off_2023}), i.e.,
\begin{itemize}
    \item \textbf{RBC:} Charge during PV surplus and discharge during load surplus such that the grid exchange becomes zero and no battery constraints are violated.
\end{itemize}

We introduce an additional timescale $\Delta c$ that denotes the control frequency at which the battery system can be operated. We set $\Delta c = \Delta gt = \Delta NI$ and with that assume that RBC can act in real-time to account for ground-truth changes. With that, we neglect phenomena like delayed battery response times observable in real-world \cite{beichter_beop_2025}.


\subsection{Two-Stage Model Predictive Control}

MPC is a proactive control strategy that relies on forecasts to determine a control action. MPC for residential battery scheduling thus typically considers upcoming prices as well as forecasts $\hat p_L$ in load and PV generation.
At each planning step, an optimization problem is solved for a selected prediction horizon. However, only the first control action is applied before a re-optimization occurs after a predefined time $\Delta s$ has passed. The optimization problem used in this study is the following:

\setlength{\jot}{2pt} 

\begin{subequations}\label{eq:opt}
\begin{alignat}{2}
\min_{p_B}\quad  C^{{imp}}&(p_G^{imp})+C^{{exp}}(p_G^{exp})+C^{{deg}}(p_B) \label{eq:obj}\\
\text{s.t. for all } k \in \mathcal{K} \nonumber \\
p_G(k)            & = \hat p_L(k) - p_B(k) \label{op:power-balance}\\
e(k+1) &= e(k) -  \Delta s \Big ( p_B^{ch}(k) \eta^{ch} + \frac{p_B^{dis}(k)}{\eta^{dis}} \Big ) \label{op:soe}\\
e^{\min}          & \le e(k) \le e^{\max} \label{op:esoebounds}\\
p_B^{\min}        & \le p_B(k) \le p_B^{\max} \label{op:pbounds}\\[4pt]
p_G(k)            & = p_G^{{imp}}(k) + p_G^{{exp}}(k) \label{op:pg_split1}\\
p_G^{imp}(k)\,p_G^{{exp}}(k) & \ge -10^{-8} \label{op:pg_split2}\\
p_G^{imp}(k) & \ge 0,\qquad p_G^{{exp}}(k) \le 0 \label{op:pg_split3}\\[4pt]
p_B(k)            & = p_B^{{ch}}(k) + p_B^{{dis}}(k) \label{op:pb_split1}\\
p_B^{ch}(k)\,p_B^{{dis}}(k) & \ge -10^{-8} \label{op:pb_split2}\\
p_B^{{dis}}(k) & \ge 0,\qquad \ p_B^{{ch}}(k) \le 0. \label{op:pb_split3}
\end{alignat}
\end{subequations}

The objective is to minimize total costs consisting of electricity import costs $C^{{imp}}$, electricity export revenues $C^{{exp}}$ and battery degradation costs $C^{{deg}}$.
In general, this optimization problem is relatively standard.  
For brevity, we refer to related works that derive and state all the necessary equations separately \cite{banfield_comparison_2020}. Nevertheless, we shortly introduce the variables and equations.
$p_B$ denotes the battery power and is the only decision variable and thus the only degree of freedom of the system.
Equation \eqref{op:power-balance} states that battery power $p_B$, grid power $p_G$ and forecasted net-load $\hat{p}_L$ need to be in balance. 
The SoE dynamics in Equation \eqref{op:soe} update the battery state $e$ over each step of length $\Delta s$, with charging/discharging efficiencies $\eta^{ch}$/$\eta^{dis}$.
The battery's physical constraints are defined by minimum/maximum SoE $e^{\min}$/$e^{\max}$ and battery power $p_B^{\min}$/$p_B^{\max}$.
Equations \eqref{op:pg_split1}-\eqref{op:pg_split3} and \eqref{op:pb_split1}-\eqref{op:pb_split3} are complementary constraints than enable the modeling of grid power as imports $p_G^{imp}$ and exports $p_G^{exp}$, as well as the split of battery power into charging/discharging power $p_B^{ch}$/$p_B^{dis}$ respectively.

All powers indexed with k represent averages over the k-th timestep. 
To formalize this, we introduce a fourth and final timescale $\Delta s$. The scheduler step length $\Delta s$ describes the time between subsequent optimization steps, commonly 60/30/15 minutes. 
In this work, $\Delta s$ may differ from the control/ground-truth/settlement timescale, so $\Delta s\neq \Delta c=\Delta gt = \Delta NI$ can hold. In those instances, MPC issues one setpoint per $\Delta s$ and cannot react to fluctuations that occur within the interval.
To address this issue, we mirror real-world implementations and add a fast-acting RBC layer on top of the MPC as a second layer to deal with intra-interval fluctuations. Two variants are investigated:
\begin{itemize}
    \item \textbf{MPC Const-Grid:} Adjust the battery power such that the grid exchange follows the MPC-planned value. The fluctuations are compensated by the battery system and the grid exchange remains constant (if possible).
    \item  \textbf{MPC Const-Bat:} Keep the battery power constant at the MPC setpoint throughout $\Delta s$, shifting all fluctuations to the grid. 
\end{itemize}

\section{Experimental Setup}
\label{sec:exp-setup}

This section contains the setup of the performed simulations. \autoref{sec:data} describes the data we use for simulation and forecasting. \autoref{sec:forecasts} contains the forecasts needed for the MPC approaches. The battery specifications and degradation costs are described in \autoref{sec:battery-specs}. \autoref{sec:prices} describes the selected import/export prices. Finally, \autoref{sec:config} describes further configurations of the experiments.

\subsection{Data}
\label{sec:data}

We build our investigation on the data presented in \cite{schlemminger_dataset_2022}, which provides high-resolution electricity consumption for several Single-Family Houses (SFH) in an enclosed district in northern Germany. We select a total of 15 buildings according to their data availability for this study. All of the buildings feature a respective heat pump. As no individual on-site PV measurements are available, we construct realistic net-load series by combining the building's load measurements with PV measurements from three nearby PV sites included in the dataset (east-, south-, and west-facing). We randomly assign one of the three orientations (and with that one of the three PV measurements) to each building.
We furthermore scale each newly created residential PV system according to its building size. For that, we map a given number of PV panels to each building according to its documented building area. The total number of PV panels varies between 16 and 38 panels per building. With that, we create a high-resolution dataset based on real-world measurements for 15 single-family buildings spanning approximately 2 years and 5 months of residential net-load data in a 1-minute resolution.

\subsection{Forecasting}
\label{sec:forecasts}
Our MPCs rely on forecasts of the upcoming net-load. 
To separate forecasting quality from timescale effects, we evaluate two forecast types:
\begin{itemize}
    \item \textbf{Ideal Forecasts:} Ideal interval-average forecasts of net-load built from ground-truth.\footnote{\textit{Ideal Forecasts} refer to the exact 60/30/15-minute averages over the next 24 hours. They specifically do not include any information about intra-interval fluctuations within each average.} 
    \item \textbf{Real Forecasts:} Data-driven realistic forecasts of net-load using Kolmogorov-Arnold Networks (KAN) \cite{liu_kan_2025}.
\end{itemize}
We predict 24 hours of residential net-load in three resolutions (60, 30, and 15 minutes), yielding a horizon of 24, 48 and 96 interval averages, respectively. 
To obtain the realistic forecasts, we employ KAN in an autoregressive manner, relying purely on historical net-load. No exogenous variables like weather forecasts or calendric features are included. We split the data chronologically and perform time-series cross-validation \cite{hyndman2018forecasting}: the first year is used for training, the second for validation, followed by a five-month test period (mid-July to mid-December). For each building, we train locally \cite{montero-manso_principles_2021}, evaluating 125 randomly sampled hyperparameter configurations from a predefined grid and select the model with the best validation performance. We then generate rolling 24-hour forecasts over the test set. 

\subsection{Battery Specifications}
\label{sec:battery-specs}

The selected battery specifications are listed in \autoref{tab:battery-params}. To obtain realistic battery degradation costs, we refer to the warranty guarantee of a specific battery model \cite{byd-box}, ensuring that the capacity retention will be above 60\% usable energy until a minimum discharged energy of 42.69 MWh is reached. We assume that the battery can be resold after 60\% maximum capacity remains for 10\% of its original value.
The degradation costs are calculated according to
\begin{align}
    c^{deg} &= \frac{\text{Investment Costs} - \text{Salvage Return}}{\text{Lifetime Discharged Energy}} \nonumber \\ 
    &= \frac{4000\text{€} - 400\text{€}}{42690\,\text{kWh}} = 0.084\text{€/kWh}.
\end{align}
Note that we do not consider a continuous degradation, meaning that during our simulation the maximum battery SoE $e^{\text{max}}$ is constant over time. This is justifiable since our evaluation is based on 5 months of simulation and during those the maximum discharged total energy is below 1400 kWh. 

We compute the battery degradation costs according to
\begin{equation}
    C^{{deg}} = c^{{deg}}\sum_{k \in \mathcal{K}} \tfrac{p_B^{dis}(k)}{\eta^{dis}}.
\end{equation}
Note that the degradation costs are derived from discharged battery power only, 
since the selected warranty defines its guarantees over the batteries total energy output, i.e., its discharged energy only.

\begin{table}[ht] 
\centering 
\caption{Selected battery specifications.} 
\label{tab:battery-params} 
\begin{tabular}{c c c c c } 
\toprule 
$e^{\min}$ [kWh] & $e^{\max}$ [kWh] & $p_{B}^{\min}$ [kW] & $p_{B}^{\max}$ [kW] & $\eta^{ch/dis}$ [\%] \\ 
\midrule 
0.0 & 13.8 & -5.0 & 5.0 & 98 \\ 
\bottomrule 
\end{tabular} 
\end{table}

\subsection{Pricing Scheme}
\label{sec:prices}

\begin{figure}[ht]
    \centering
    \includegraphics[width=0.95\linewidth]{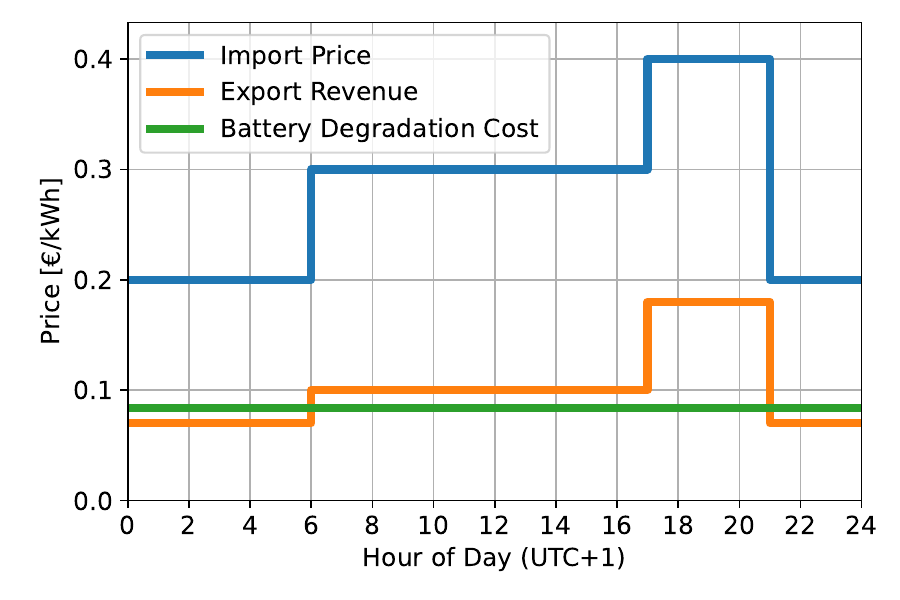}
    \caption{Selected TOU prices for imports $c^{imp}$ and exports $c^{exp}$ with constant battery degradation costs $c^{{deg}}$.}
    \label{fig:prices-tou}
\end{figure}

We calculate the total import/export costs according to
\begin{align}
    C^{imp} &= \sum_{k \in \mathcal{K}} c^{imp}(k) p_G^{imp}(k) \\
    C^{exp} &= \sum_{k \in \mathcal{K}} c^{exp}(k) p_G^{exp}(k), 
\end{align}
with $c^{imp}(k)$ and $c^{exp}(k)$ representing the €/kWh prices for imports and exports respectively.
\autoref{fig:prices-tou} depicts the selected tariff structure for the grid exchange and battery degradation. 
We orient the exchange prices at a typical three-level residential Time-of-Use (TOU) pricing scheme, characterized by lower prices during nighttime and a high-pricing period during the typical evening high net-load peak.
Generally, the more prices vary throughout the day, the more an MPC approach can leverage them to minimize costs. Furthermore, to align the prices with EUs directive of promoting self-sufficiency \cite{varga2024changing}, we select a tariff structure that specifically does not incentivize arbitrage. That is, $c^{imp}(t_0) > c^{exp}(t_1) \forall t_0, t_1$.
Finally, it can be seen that at night, discharging the battery system to export energy is a net loss for the residents because the degradation costs surpass the export revenue.

\subsection{Configuration}
\label{sec:config}

We simulate battery control for five consecutive months of data from mid-July to mid-December 2020 for 15 buildings. 
We design two main experiments. 
For each experiment, we simulate battery control using the different battery control models described in \autoref{sec:RBC-MPC}. Each model is furthermore evaluated three times using different scheduling interval lengths:
\begin{equation}
    \Delta s \in \{60, 30, 15\}\ \text{min},
\end{equation}
with the respective 60/30/15 min forecasts described in \autoref{sec:forecasts}.
The difference between the two main experiments is the ground-truth data resolution $\Delta gt$:
\begin{itemize}
    \item \textbf{Fully Averaged Evaluation}: The ground-truth data resolution equals the scheduling interval length, i.e., $\Delta gt = \Delta s$. 
    Only coarser averages with length $\Delta s$ are used to assess the models' performances.
    This reproduces a typical evaluation strategy as described in \autoref{sec:related-works}.
    \item \textbf{Fine-Resolution Evaluation}: The ground-truth data resolution is fixed at $\Delta gt = 1 \text{min}$. 
    The models' performances are based on the aggregation of 1-min averages. 
    All models now rely on faster-acting rule-based behavior with control frequency $\Delta c = \Delta gt \neq \Delta s$.
    This evaluation strategy aims to obtain more realistic performance assessments under Net Billing settlement rules.\footnote{Note that in reality, Net Billing stems from instantaneous measurements and not from 1 minute averages. To be precise, with this setup, we aim to approximate Net Billing by investigating the performance under Net Metering with a Netting Interval $\Delta NI = 1 \text{min}$.}
\end{itemize}

All models and evaluations are solved in Python using the optimization modeling language Pyomo \cite{bynum_pyomo_2021} with IPOPT \cite{wachter2006implementation} as the solver. The implementation is published as an open-access repository along with this paper. All presented results can be reproduced using the provided code.\footnote{Code will be published as an open-access repository under CC-BY license upon acceptance.}

\section{Results}
\label{sec:results}

The experiment results are reported in \autoref{tab:res-gt-s} and \autoref{tab:res-gt-1}. Both tables contain mean results across the 15 buildings. \autoref{tab:res-gt-s} corresponds to the Fully Averaged Evaluation with $\Delta s = \Delta gt$. 
Here, the terms "(Ideal)" or "(Real)" indicate whether the results are achieved using ideal or realistic forecasts. Additionally, the label "MPC (Ideal)" does not contain a specified MPC rule (neither Const-Grid nor Const-Bat), because both rules yield the same result. \autoref{tab:res-gt-1} shows the Fine-Resolution Evaluation with $\Delta s \neq \Delta gt$.

Each table reports imported energy, exported energy, their associated costs, discharged energy and degradation costs. Total Costs is the sum of the electricity bill and degradation costs and serves as the main performance indicator for the present work. The relative performance (Rel. Perf.) describes the change of total costs relative to RBC. Finally, Ranking denotes the average placement of the respective model over 15 buildings.\footnote{For example, a model with a Ranking of 1.00 means that for 15 out of 15 buildings, this model performed better than the others.} In the following, the results of the tables are visualized and discussed.

\begin{table*}[ht]
\centering
\caption{\textbf{Fully Averaged Evaluation:} \boldmath $\Delta s = \Delta gt$}
\label{tab:res-gt-s}
\scriptsize
\setlength{\tabcolsep}{6pt}
\begin{tabular}{l c c c c c c c c c c}
\toprule
\textbf{Model} & \textbf{Imp.} & \textbf{Costs Imp.} & \textbf{Exp.} 
               & \textbf{Costs Exp.} & \textbf{E Dis.} & \textbf{Degrad.}
               & \textbf{Bill} & \textbf{Total Costs} & \textbf{Rel. Perf.} & \textbf{Ranking} \\
              & [kWh] & [€] & [kWh] & [€] & [kWh] & [€] & [€] & [€] & [\%] & [-] \\
\midrule
\midrule
\multicolumn{10}{c}{\textbf{Schedule Step \boldmath$\Delta s=$ 60min \quad Ground-Truth Resolution \boldmath $\Delta gt=$ 60min}} \\
\midrule
MPC (Ideal) & 1191.19 & \textbf{250.77} & 1149.04 & 123.78 & 1039.20 & 87.29 & \textbf{125.01} & \textbf{212.31} & -12.55 & 1.00 \\
RBC         & \textbf{1132.74} & 322.04 & 1122.48 & 120.93 & \textbf{496.66} & \textbf{41.72} & 201.05 & 242.77 & - & 2.20 \\
MPC (Real) Const-Grid & 1304.44 & 298.72 & \textbf{1262.46} & \textbf{138.74} & 1078.74 & 90.61 & 158.35 & 248.96 & 2.55 & 2.80 \\
\midrule
\midrule
\multicolumn{10}{c}{\textbf{Schedule Step \boldmath$\Delta s=$ 30min \quad Ground-Truth Resolution \boldmath $\Delta gt=$ 30min}} \\
\midrule
MPC (Ideal) & 1192.76 & \textbf{251.40} & 1149.64 & 124.13 & 1065.63 & 89.51 & \textbf{125.31} & \textbf{214.83} & -12.41 & 1.00 \\
RBC & \textbf{1135.01} & 322.61 & 1123.68 & 121.23 & \textbf{523.20} & \textbf{43.95} & 201.32 & 245.27 & - & 2.20 \\
MPC (Real) Const-Grid & 1291.19 & 295.27 & \textbf{1247.55} & \textbf{135.98} & 1128.68 & 94.81 & 157.72 & 252.53 & 2.96 & 2.80 \\
\midrule
\midrule
\multicolumn{10}{c}{\textbf{Schedule Step \boldmath$\Delta s=$ 15min \quad Ground-Truth Resolution \boldmath $\Delta gt=$ 15min}} \\
\midrule
MPC (Ideal) & 1193.69 & \textbf{251.88} & 1149.59 & 124.29 & 1090.02 & 91.56 & \textbf{125.63} & \textbf{217.19} & -12.30 & 1.00 \\
RBC & \textbf{1136.74} & 323.08 & 1124.43 & 121.37 & \textbf{547.49} & \textbf{45.99} & 201.65 & 247.64 & - & 2.13 \\
MPC (Real) Const-Grid & 1289.26 & 294.92 & \textbf{1243.67} & \textbf{135.10} & 1183.89 & 99.45 & 158.30 & 257.75 & 4.08 & 2.87 \\
\bottomrule
\end{tabular}
\end{table*}

\begin{table*}[ht]
\centering
\caption{\textbf{Fine-Resolution Evaluation:} \boldmath $\Delta s \neq \Delta gt = 1 \textbf{min}$.}
\label{tab:res-gt-1}
\scriptsize
\setlength{\tabcolsep}{6pt}
\begin{tabular}{l c c c c c c c c c c}
\toprule
\textbf{Model} & \textbf{Imp.} & \textbf{Costs Imp.} & \textbf{Exp.} 
               & \textbf{Costs Exp.} & \textbf{E Dis.} & \textbf{Degrad.}
               & \textbf{Bill} & \textbf{Total Costs} & \textbf{Rel. Perf.} & \textbf{Ranking} \\
              & [kWh] & [€] & [kWh] & [€] & [kWh] & [€] & [€] & [€] & [\%] & [-] \\
\midrule
\midrule
\multicolumn{10}{c}{\textbf{Schedule Step \boldmath$\Delta s=$ 60min --------- Ground-Truth Resolution \boldmath $\Delta gt=$ 1min}} \\
\midrule
MPC (Ideal) Const-Grid & 1207.75 & \textbf{256.79} & 1153.98 & 124.28 & 1329.11 & 111.65 & \textbf{130.55} & \textbf{242.19} & -3.76 & 1.20 \\
RBC & \textbf{1139.61} & 323.93 & 1125.71 & 121.54 & \textbf{587.12} & \textbf{49.32} & 202.33 & 251.65 & - & 1.80 \\
MPC (Real) Const-Grid & 1318.83 & 303.31 & 1266.00 & 139.08 & 1347.95 & 113.23 & 162.60 & 275.83 & 9.61 & 2.00 \\
MPC (Ideal) Const-Bat & 1671.39 & 404.92 & \textbf{1629.33} & \textbf{181.42} & 1039.26 & 87.30 & 221.54 & 308.83 & 22.72 & 4.00 \\
\midrule
\midrule
\multicolumn{10}{c}{\textbf{Schedule Step \boldmath$\Delta s=$ 30min --------- Ground-Truth Resolution \boldmath $\Delta gt=$ 1min}} \\
\midrule
MPC (Ideal) Const-Grid & 1202.98 & \textbf{255.36} & 1150.83 & 124.25 & 1288.74 & 108.25 & \textbf{129.14} & \textbf{237.40} & -5.66 & 1.13 \\
RBC & \textbf{1139.61} & 323.93 & 1125.71 & 121.54 & \textbf{587.12} & \textbf{49.32} & 202.33 & 251.65 & - & 1.87 \\
MPC (Real) Const-Grid & 1300.82 & 298.40 & 1248.83 & 136.09 & 1336.25 & 112.25 & 160.75 & 272.99 & 8.48 & 3.07 \\
MPC (Ideal) Const-Bat & 1583.96 & 377.07 & \textbf{1540.85} & \textbf{171.06} & 1065.41 & 89.49 & 204.05 & 293.54 & 16.65 & 3.93 \\
\midrule
\midrule
\multicolumn{10}{c}{\textbf{Schedule Step \boldmath$\Delta s=$ 15min --------- Ground-Truth Resolution \boldmath $\Delta gt=$ 1min}} \\
\midrule
MPC (Ideal) Const-Grid & 1199.80 & \textbf{254.38} & 1149.56 & 124.28 & 1242.28 & 104.35 & \textbf{128.14} & \textbf{232.49} & -7.61 & 1.00 \\
RBC & \textbf{1139.61} & 323.93 & 1125.71 & 121.54 & \textbf{587.12} & \textbf{49.32} & 202.33 & 251.65 & - & 2.13 \\
MPC (Real) Const-Grid & 1294.77 & 296.78 & 1243.67 & 135.08 & 1321.32 & 110.99 & 160.20 & 271.19 & 7.76 & 3.27 \\
MPC (Ideal) Const-Bat & 1472.98 & 341.53 & \textbf{1428.90} & \textbf{157.74} & 1089.99 & 91.56 & 181.84 & 273.40 & 8.64 & 3.60 \\
\bottomrule
\end{tabular}
\end{table*}

We begin with a few general observations. 
At first, we want to draw attention to the Fully Averaged Evaluation in \autoref{tab:res-gt-s}. Within each configuration, models are ordered from lowest to highest total cost. For each configuration $\Delta s$, MPC (Ideal) ranks best, RBC second, and MPC (Real) third.  
An observation that may seem counterintuitive is that total costs increase for all models as the scheduling step $\Delta s$ decreases. Here, $\Delta gt = \Delta s$, so a smaller $\Delta s$ also leads to a finer ground-truth that reveals more fluctuations that need to be handled by the battery. This increases costs, because $c^{imp}>c^{exp}$ holds.
Further, the relative performance of MPC (Ideal) over RBC is -12.41\% for $\Delta s=30\text{min}$, aligning with the -12.9\% cost reduction reported in \cite{banfield_comparison_2020}.\footnote{To be precise, \cite{banfield_comparison_2020} only reports absolute values, i.e., 759\$ for RBC and 661\$ for MPC (Ideal), which yields a reduction of total costs by -12.9\%.}

In the Fine-Resolution Evaluation in \autoref{tab:res-gt-1}, the ranking of the models also is consistent, with an MPC (Ideal) performing best. Here, we want to highlight two observations: First, RBC's total costs of 251.65€ is constant across $\Delta s$. RBC does not schedule, it reacts to the ground-truth, which is fixed at $\Delta gt = 1\text{min}$ here. Second, MPC (Ideal) Const-Bat performs worst. Fixing the battery setpoint over each scheduling interval shifts all intra-interval fluctuations to the grid and induces frequent import/export swings. The results highlight how important intra-interval fluctuations are under Net Billing. Even with imperfect forecasts, combining MPC with an effective fast rule that compensates for fluctuations can beat an ideal-forecast MPC paired with a poor fast rule.

\subsection{Impact of Battery Control Strategies and Imperfect Forecasting on Cost Composition}

\begin{figure}[ht]
    \centering
    \includegraphics[width=0.98\linewidth]{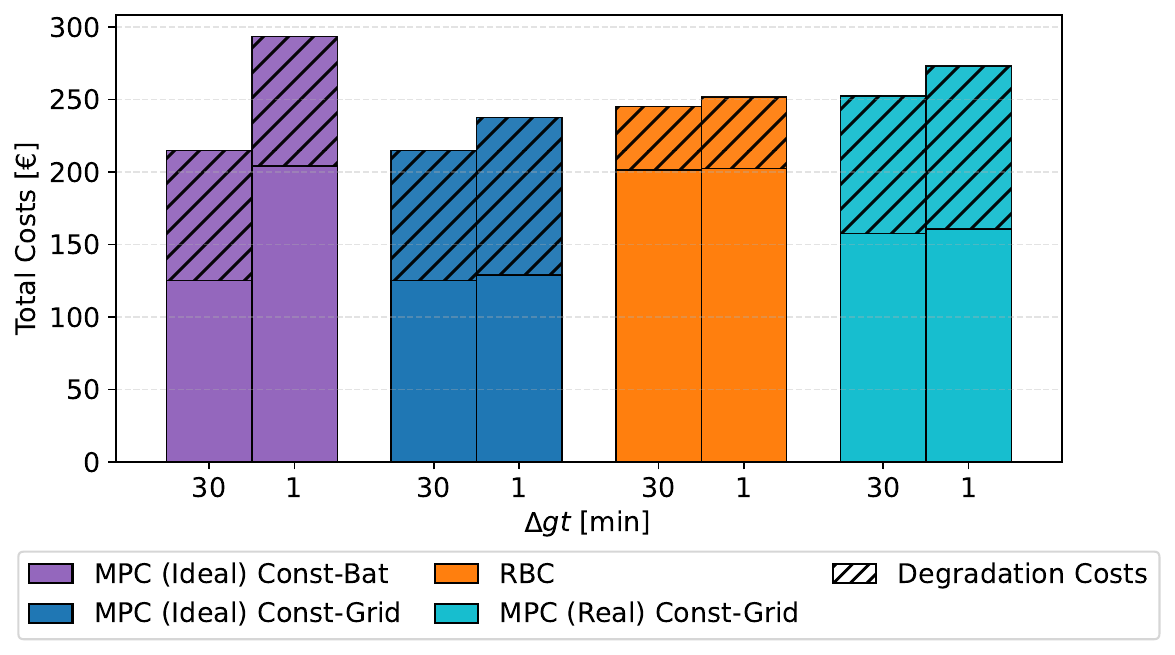}
    \caption{Total costs for two different ground-truth resolutions $\Delta gt$ with an MPC scheduler step length of $\Delta s = 30\text{min}$. All MPC approaches lead to comparatively high battery degradation costs when compared to RBC.}
    \label{fig:cost-composition}
\end{figure}

In \autoref{fig:cost-composition}, the total costs of the different models for a scheduling step length of $\Delta s=30\text{min}$ are displayed. 
All MPC schedulers contain a fair amount of degradation costs, whereas the majority of the costs of RBC stem from electricity costs. 
The main reason for that is that MPCs actively use the battery system to leverage prices.
A battery controlled via RBC on the other hand is in many cases fully charged during summer, and fully discharged during winter. This reduces battery usage at the cost of higher electricity bills.

Furthermore, all approaches lead to higher costs at the finer time resolution $\Delta gt=1$min. These additional costs mainly stem from battery degradation, except for the constant-battery MPC. 
The degradation costs of MPC Const-Grid with ideal and real forecasts are similar, but electricity costs differ.
Finally, while perfect-forecast MPC can perform better than RBC, the margin is relatively small.
With realistic forecasts, MPC can underperform RBC in our setup, making the practical advantage of MPC questionable for this scenario.


\subsection{Systematic Underestimation of Total Costs under Averages}

\begin{figure}[ht]
     \centering
     \includegraphics[width=1.0\linewidth]{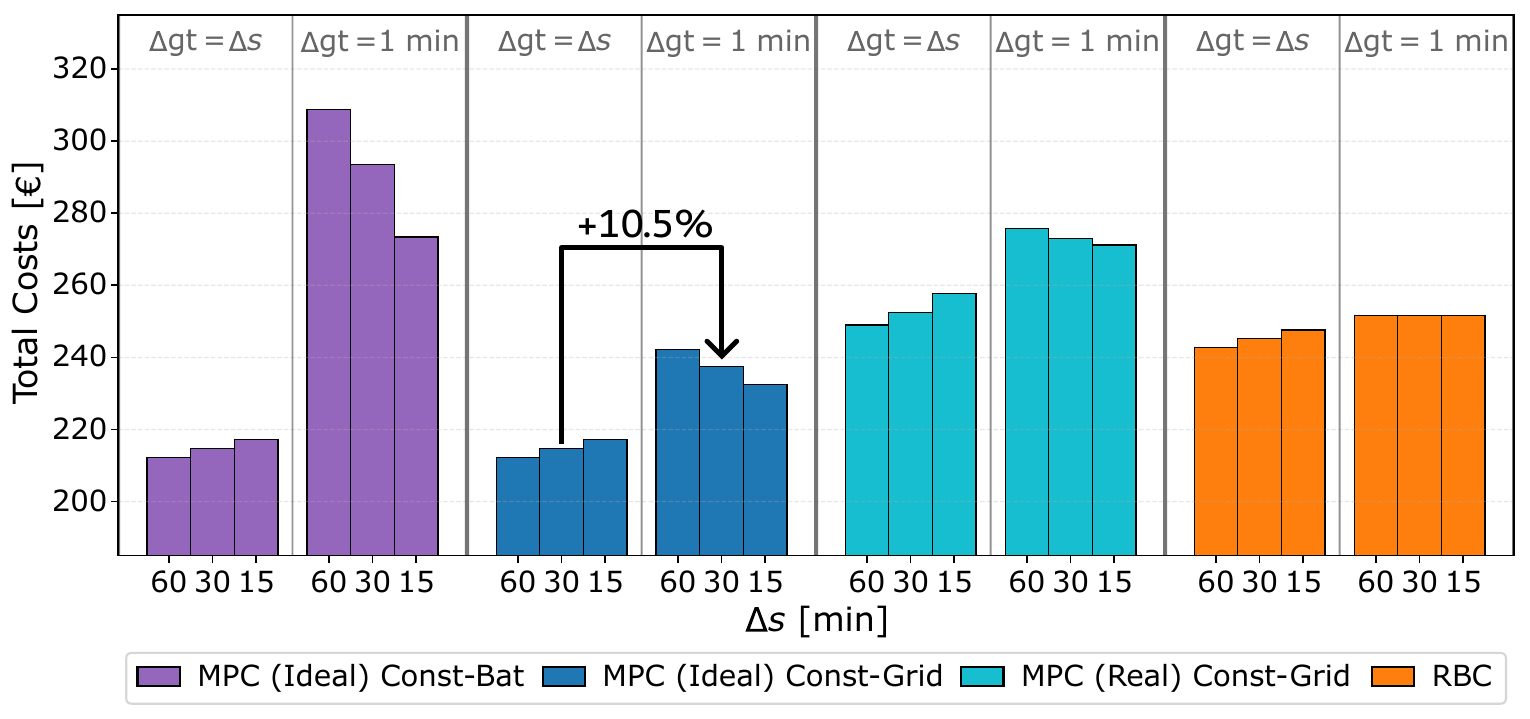}
     \caption{Total costs for varying scheduler step lengths $\Delta s$. Left-hand blocks of a specific color indicate Fully Averaged Evaluations ($\Delta gt = \Delta s$). Right-hand blocks indicate Fine-Resolution Evaluations ($\Delta gt = 1\text{min}$). Averaging hides costs. Throughout all performed experiments, average evaluations using $\Delta gt = \Delta s$ underestimate costs. The bigger $\Delta s$, the bigger the underestimation.}
     \label{fig:total-costs-all}
 \end{figure}

\autoref{fig:total-costs-all} shows the total costs for all configurations. 
Here, it can be seen that averaging leads to cost underestimation: in each case, the bars in the left block ($\Delta gt = \Delta s$) are lower than their counterparts on the right.
See for example MPC (Ideal) Const-Grid for $\Delta s = 30\text{min}$: Averaging over 30-min intervals underestimates total costs by $(237.40\text{€} - 214.83\text{€})/214.83\text{€} = 10.5\%$ as shown in \autoref{fig:total-costs-all}. 
\autoref{tab:relative-underestimation} contains the relative underestimation of total costs for all configurations. It can be seen that the underestimations of costs for RBC are substantially lower than the ones of MPC. This effect is further analyzed in the following subsection.

\begin{table}[ht]
  \centering
  \caption{Relative underestimation of Total Costs. By how much are costs underestimated when evaluations are performed using $\Delta gt = \Delta s$ instead of $\Delta gt = 1 \text{min}$.}
  \label{tab:relative-underestimation}
  \setlength{\tabcolsep}{4pt}
  \renewcommand{\arraystretch}{1.15}
  \begin{tabular}{
      p{0.19\linewidth}  
      >{\centering\arraybackslash}p{0.17\linewidth}
      >{\centering\arraybackslash}p{0.17\linewidth}
      >{\centering\arraybackslash}m{0.10\linewidth} 
      >{\centering\arraybackslash}p{0.17\linewidth}
    }
    \toprule
    & MPC (Ideal) Const-Bat & MPC (Ideal) Const-Grid & RBC & MPC (Real) Const-Grid \\
    \midrule
    $\Delta s = 60\,\mathrm{min}$ & 45.5\% & 14.1\% & 3.7\% & 10.8\% \\
    $\Delta s = 30\,\mathrm{min}$ & 36.6\% & 10.5\% & 2.6\% & 8.1\% \\
    $\Delta s = 15\,\mathrm{min}$ & 25.9\% & 7.0\%  & 1.6\% & 5.2\% \\
    \bottomrule
  \end{tabular}
\end{table}

\subsection{Systematic Bias in MPC vs. RBC Comparison under Averages}

To omit forecasting-related effects and because the constant battery setpoint rule underperforms, we focus in this subsection exclusively on the comparison of MPC (Ideal) Const-Grid versus RBC. 

When costs are computed from averages (\autoref{tab:res-gt-s}), MPC performs clearly better. With a scheduling step length of $\Delta s=30\text{min}$, the averaged evaluation shows a 30.44€ reduction of costs relative to RBC. Under Net Billing, however, real operation includes intra-interval fluctuations that cause costs. When the same controllers are evaluated at finer resolution (\autoref{tab:res-gt-1}), MPC's advantage drops to 14.25€. This is a shrinkage of MPC's presumed advantage of $1-14.25/30.44 = 53\%$. For brevity, we call this metric the shrinkage. \autoref{tab:shrinkages-average} contains the shrinkages for all three scheduling step lengths $\Delta s$. The potential advantage of MPC shrinks by as much as 69\%, which is far from negligible. If decisions are based on averaged evaluations alone, a reduction of this magnitude could change the preferred control strategy.

\begin{figure}[ht]
    \centering
    \includegraphics[width=0.95\linewidth]{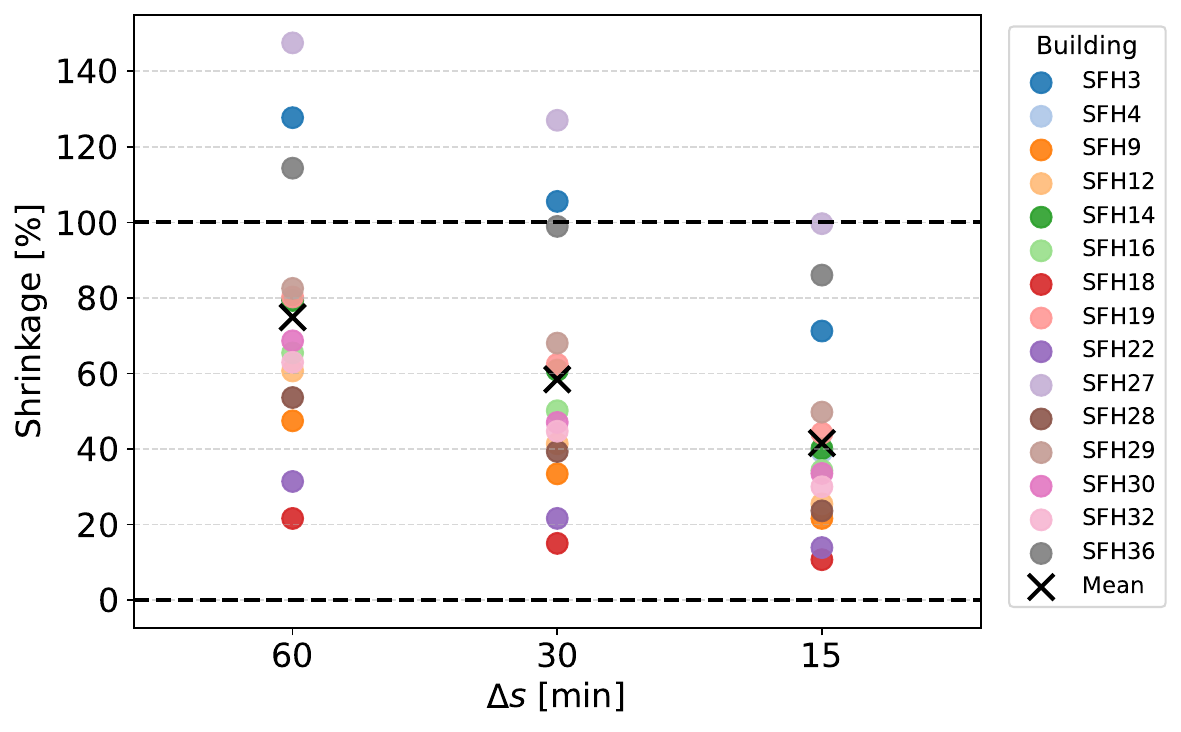}
    \caption{Individual building shrinkage for different scheduling step lengths $\Delta s$. The higher the shrinkage, the more of MPCs presumed advantage over RBC decreases once fluctuations are considered. Shrinkages over 100\% indicate that all advantage is lost and RBC leads to better performance than an MPC with ideal forecasts.}
    \label{fig:shrinkage-individual-b}
\end{figure}

Up until now, all shown results are based on the average over the 15 simulated buildings. \autoref{fig:shrinkage-individual-b} shows the shrinkage for each building individually.\footnote{Note that the mean shrinkage visualized in \autoref{fig:shrinkage-individual-b} does not perfectly equal the mean shrinkage in \autoref{tab:shrinkages-average}. They differ because \autoref{tab:shrinkages-average} contains the shrinkage of an average building, whereas the figure plots the average shrinkage of a building.} Building-level variability is high, and the general trend persists: as $\Delta s$ decreases, shrinkage becomes smaller.
However, at the same time, forecasting at shorter steps is more challenging, making the assumption of ideal forecasts less justifiable.

\begin{table}[ht] 
\centering 
\caption{Shrinkage---How much of MPC's advantage over RBC is lost once fluctuations are considered.} 
\label{tab:shrinkages-average} 
\begin{tabular}{c c c c} 
\toprule 
$\Delta s$ & 60 min & 30 min & 15 min \\ 
\midrule 
Shrinkage & 69\% & 53\% & 37\% \\ 
\bottomrule 
\end{tabular} 
\end{table}

Moreover, in several cases the shrinkage exceeds 100\%, i.e., the expected advantage is not only decreased but reverses, and RBC outperforms MPC even under perfect forecasts for individual buildings.
For instance, for SFH3 at $\Delta s = 60\text{min}$, the Fully Averaged Evaluation indicates MPC savings of 18\% (30€) relative to RBC. However, the Fine-Resolution Evaluation shows RBC is better by 4\% (8€). Hence, the presumed advantage shrinks by 128\%. In this specific case, an average-based evaluation would recommend the wrong controller. 

Additionally, we want to highlight that not a single building has a shrinkage below 0\%. 
Across all 45 simulations (15 buildings × 3 scheduler step lengths $\Delta s$), the advantage of MPC never increases relative to RBC when a more realistic evaluation is applied.
The underlying reason for that is that intra-interval fluctuations under RBC are more often absorbed by the grid. MPC schedules typically exhibit smoother behavior throughout the day, reaching full charge or depletion at later times, and thus use the battery system more to counteract fluctuations, leading to increased degradation costs.
To conclude, coarse-average evaluations overlook fluctuation effects and thereby systematically favor MPC over RBC. Identifying and quantifying this bias is a key contribution of the present work.

\section{Discussion}
\label{sec:discussion}

In our evaluation, Net Billing is approximated with 1-min averages. Thus, we underestimate costs even for the Fine-Resolution Evaluation and with that also understate the reported effects and drawn conclusions. Under this approximation, our results show that a fair comparison of battery scheduling methods must consider high-resolution fluctuations and the controller’s response to them. By contrast, studies focusing on Net Metering settlements or on tracking a committed schedule can reasonably use a simpler coarse-average evaluation that ignores within-interval variability.

\textbf{Limitations:}
Degradation is modeled per discharged kWh only. We do not consider SoE-dependent, calendar or rate-dependent aging effects. High SoE is known for accelerating aging, which is why some scheduling approaches specifically try to delay charging \cite{bui_study_2021}. The negligence of this effect skews our results in favor of RBC. However, this just affects the relative performance between RBC and MPC and does not impair our main contribution, i.e., the effects of averaging on drawn conclusions.
Further, MPC has a larger share of degradation costs than RBC. If battery prices and the implied €/kWh degradation cost decrease in the future, costs for MPC would drop more than for RBC, improving the relative advantage of MPC.

\textbf{Implications for practice and policy:}
Even with time-varying prices, RBC provides a strong baseline. This is strengthened by the fact that we do not consider additional installation, integration, and maintenance costs typically associated with real-world MPC deployments. RBC is also straightforward to enhance with simple rules (e.g., pre-charging at night to a target SoE during winter). 
Yet, RBC is generally less grid-friendly than MPC: midday PV spikes, for instance, are more likely to be fully shifted to the grid. This hints at a misalignment between common Net Billing tariffs and grid objectives. If policymakers want prosumers to adopt MPC and deliver smoother exchanges, pricing and settlement should be designed so MPC realistically outperforms RBC by a clear margin. Levers include widening price variations, introducing capacity-based tariffs, and also adjusting settlement rules, as we have shown in this work. For example, regulators could employ Net Metering with a short netting interval (such as 15 minutes) instead of Net Billing.

\section{Conclusion}
\label{sec:conclusio}

This paper investigates performance of Model Predictive Control (MPC) and  Rule-Based Control (RBC) for residential battery systems regarding the minimization of costs and battery degradation. We investigate timescale effects, specifically the justifiability of performing evaluations based on 15/30/60-min averages.
We find that under Net Billing, averaged evaluations underestimate total costs and overestimate the advantages of MPC over RBC significantly. 

For a better scientific practice, we recommend researchers to specify their settlement assumptions and state clearly when, if, and why the investigated settlement rule does not match the true underlying settlement rule. To elaborate, if performance is measured on 30-min averages, the investigated results actually represent performance on a Net Metering scheme with an NI of 30 minutes, instead of Net Billing. 

To conclude, averaging hides fluctuations that cause real costs and impact different battery control methods differently. Thus, in order to obtain a solid foundation for a fair comparison between different battery control approaches, accounting for high frequency fluctuations is a necessity that should not be neglected by only considering coarse averages. 
Future work can focus on investigating other fast-acting rule strategies that might synergize better with slow-planning MPCs.

\section*{Acknowledgments}
The authors thank the Helmholtz Association for the support under the "Energy System Design" program and the German Research Foundation as part of the Research Training Group 2153 "Energy Status Data: Informatics Methods for its Collection, Analysis and Exploitation".

\section*{Declaration}
For this work, the authors used ChatGPT during writing and coding to improve clarity and fluency. After the usage, the authors reviewed and edited the content as needed and take full responsibility for the content of the published article. 

\bibliographystyle{elsarticle-num} 
\bibliography{references_new}

\end{document}